\documentclass[letterpaper, 10 pt, conference]{ieeeconf}%
\usepackage{graphics}
\usepackage{epsfig}
\usepackage{mathptmx}
\usepackage{amsmath}
\usepackage{amssymb}
\usepackage{setspace}
\usepackage[colorinlistoftodos]{todonotes}
\usepackage[colorlinks=true, allcolors=blue]{hyperref}
\usepackage{caption}
\usepackage{subcaption}
\usepackage{hyperref}
\usepackage{xcolor,colortbl}
\usepackage{setspace}
\usepackage{float}
\usepackage{cite}
\usepackage{graphicx}
\usepackage[space]{grffile}
\usepackage{amsfonts}%
\providecommand{\U}[1]{\protect\rule{.1in}{.1in}}
\IEEEoverridecommandlockouts
\begin{document}
\bstctlcite{MyBSTcontrol}

\title{{\LARGE \textbf{Eco-Routing of Plug-In Hybrid Electric Vehicles in
Transportation Networks }}}
\author{Arian Houshmand$^{1}$ and Christos G. Cassandras$^{1}$\thanks{*This work was
supported in part by NSF under grants ECCS-1509084, IIP-1430145, and
CNS-1645681, by AFOSR under grant FA9550-12-1-0113, by ARPA-E’s NEXTCAR program under grant DE-AR0000796, and by Bosch and MathWorks.}\thanks{$^{1}$The authors are with
Division of Systems Engineering, Boston University, Brookline, MA 02446 USA
\texttt{{\small arianh@bu.edu; cgc@bu.edu}}}}
\maketitle

\begin{abstract}
We study the problem of eco-routing Plug-In Hybrid Electric Vehicles (PHEVs)
to minimize the overall energy consumption costs. Unlike the traditional
Charge Depleting First (CDF) approaches in the literature where the power-train control strategy
is fixed, we propose a Combined Routing and Power-train Control (CRPTC)
algorithm which can simultaneously calculate the optimal energy route as well
as the optimal power-train control strategy. To validate our method, we apply
our eco-routing algorithm to a subnetwork of the Eastern Massachusetts (EMA)
transportation network using actual traffic data provided by the Boston Region
Metropolitan Planning Organization. As an alternative benchmark, we also
simulate the traffic behavior of the network using the extracted flow data
from the aforementioned traffic dataset. We show that the CRPTC approach outperforms the
traditional CDF approach and we quantify the trade-off between saving energy
and time in using eco-routing algorithms.


\end{abstract}

\thispagestyle{empty} \pagestyle{empty}



\section{INTRODUCTION}

\label{sec: Intro} Due to environmental concerns and high gas prices, there
has been an increasing interest in vehicles using alternative energy sources such as
electric vehicles (EV). However, EV adoption is limited by the all-electric
range (AER) considering the battery capacity in current EVs. In this respect,
Plug-In Hybrid Electric Vehicles (PHEVs) are viewed as a suitable alternative,
as they can overcome range limitations by using both gas and electricity.
Moreover, it is possible to decrease the energy
consumption cost and the carbon footprint of PHEVs using smart eco-routing and
power-train control strategies.

Unlike traditional vehicle routing algorithms which seek to find the minimum
time or shortest path routes
\cite{bertsekas_dynamic_1995,braekers_vehicle_2016,toth_vehicle_2002},
eco-routing algorithms seek the paths that minimize the total energy
consumption cost. Several routing algorithms have been proposed in the
literature for conventional vehicles which are capable of finding
the energy-optimal paths using historical and online traffic data
\cite{barth_environmentally-friendly_2007,boriboonsomsin_eco-routing_2012,andersen_ecotour:_2013,yao_study_2013,yang_stochastic_2014}%
. Kubicka et al \cite{kubicka_performance_2016} performed a study to compare the objective values proposed in the eco-routing literature and showed that the performance of eco-routing algorithms is highly dependent on the method
used to calculate the traveling cost of each link. Although eco-routing of
conventional vehicles is well studied, there is little research that addresses
the case of PHEVs \cite{guanetti_control_2018}. Jurik et al
\cite{cela_energy_2014} have addressed the eco-routing problem for HEVs based
on the vehicle longitudinal dynamics. Sun et al \cite{sun_save_2016} and Qiao
et al \cite{qiao_vehicle_2016} proposed the CDF approach to address
the eco-routing for PHEVs. Furthermore, in \cite{sun_save_2016}, the authors have
shown that energy-optimal paths typically take more time compared to the
fastest route.

The contributions of this paper are summarized as follows. After reviewing the traditional CDF eco-routing approach, we propose a Hybrid-LP Relaxation algorithm to solve this problem by reducing it to a Linear Programming (LP) problem which guarantees convergence to a global optimum.
Moreover, based on the energy
model definition in Section \ref{sec: energy model}, we propose a Combined
Routing and Power-train Control (CRPTC) eco-routing algorithm for PHEVs which
can simultaneously find the optimal energy route as well as the optimal
power-train control strategy for switching between charge depleting (CD) and
charge sustaining (CS) modes. Unlike the previous methods where the  
power-train (PT) control strategy was considered a priori
\cite{sun_save_2016,qiao_vehicle_2016}, we do not make such an assumption and
we let the optimizer choose the optimal control strategy. We formulate the
problem as a mixed integer linear programing (MILP) problem and use actual
traffic data from the Eastern Massachusetts transportation network (provided
by the Boston Region Metropolitan Planning Organization) to validate the
performance of our algorithm. As an alternative to such historical data, we
also use the SUMO simulator to investigate traffic outcomes. We show that the 
CRPTC approach can lead to improved energy savings compared to the CDF
approach while using the same energy models as in \cite{qiao_vehicle_2016}.
Finally to assess the performance of the CRPTC approach, we compare the energy
cost and traveling time of the energy-optimal route with the fastest route and
the routes obtained from actual traffic data. As in \cite{sun_save_2016}, we
show the trade-off between saving energy and time in Section
\ref{sec: numerical results}.

The remainder of this paper is organized as follows. The PHEV energy
consumption model is presented in Section \ref{sec: energy model}. A MILP
formulation is proposed in Section \ref{sec: single vehicle routing} to solve
the eco-routing problem. In Section \ref{sec: numerical results}, we use
actual historical data to validate the performance of the algorithm and also
include simulation results to compare the energy cost and traveling times of
different routing algorithms. Finally, conclusions and further research
directions are outlined in Section \ref{sec: conclusions}.


\section{PHEV Energy Consumption Modeling}

\label{sec: energy model} 
Unlike conventional vehicles where
it is possible to use the empirical equations based on velocity and
acceleration of the vehicle to estimate fuel consumption cost
\cite{kamal_model_2013}, estimating a PHEV's fuel consumption is a more involved process which over a finite time horizon can be expressed as follows:
\begin{equation}
\int_{t_{0}}^{t_{f}}(C_{gas}\dot{m}_{gas}(t)+C_{ele}P_{batt}(t))dt
\label{eqn:cost function}%
\end{equation}
where $\dot{m}_{gas}$ is the fuel consumption rate, and $P_{batt}$ is the
total electrical power used/generated by the motor/generator units. Moreover,
$C_{gas}$ ($\$/gallon$) and $C_{ele}$ ($\$/kWh$) are the cost of gas and
electricity, respectively. 

Due to the nature of our problem and to avoid
unnecessary complexities, we use a simplified model proposed by Qiao et al \cite{qiao_vehicle_2016} to calculate $\dot{m}_{gas}$ and $P_{batt}$ in our eco-routing problem formulation.
Instead of using real time driving data for a
targeted vehicle, they calculate the average $\dot{m}_{gas}$ and $P_{batt}$
per mile for different drive cycles using the software package \textit{Autonomie/PSAT}. In this method, they consider two driving modes for a
PHEV: charge-depleting (CD) and charge sustaining (CS). The CD mode refers to
the phase where the PHEV acts like an EV and consumes all of its propulsion
energy from the battery pack. Once the state of the charge (SOC) of the
battery reaches a target value, it switches to the CS mode in which the
vehicle starts using the internal combustion engine as the main propulsion
system and the battery and electric motors are only used to improve fuel
economy as in HEVs \cite{karabasoglu_influence_2013}.

\begin{figure}[h]
\centering
\begin{subfigure}{.19\textwidth}
\centering
\includegraphics[width=1\linewidth]{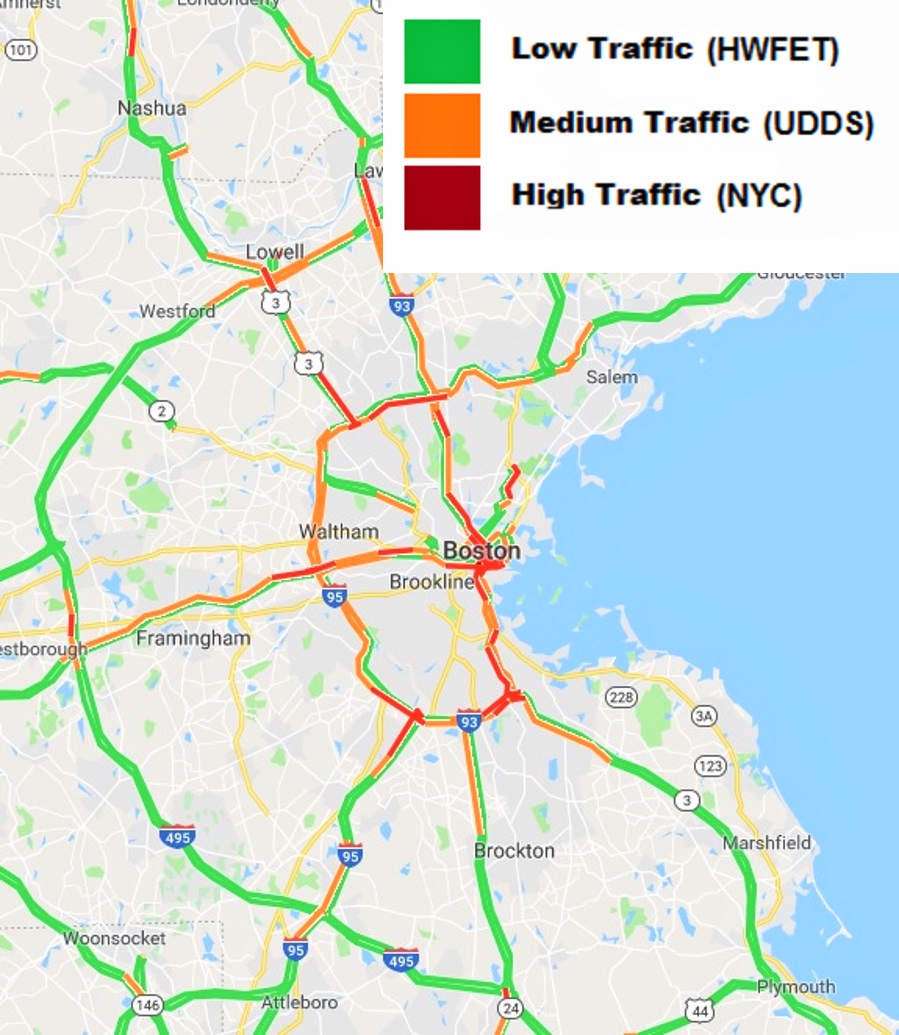}
\caption{}
\label{fig:EMA-Google}
\end{subfigure}\begin{subfigure}{.19\textwidth}
\centering
\includegraphics[width=01\linewidth]{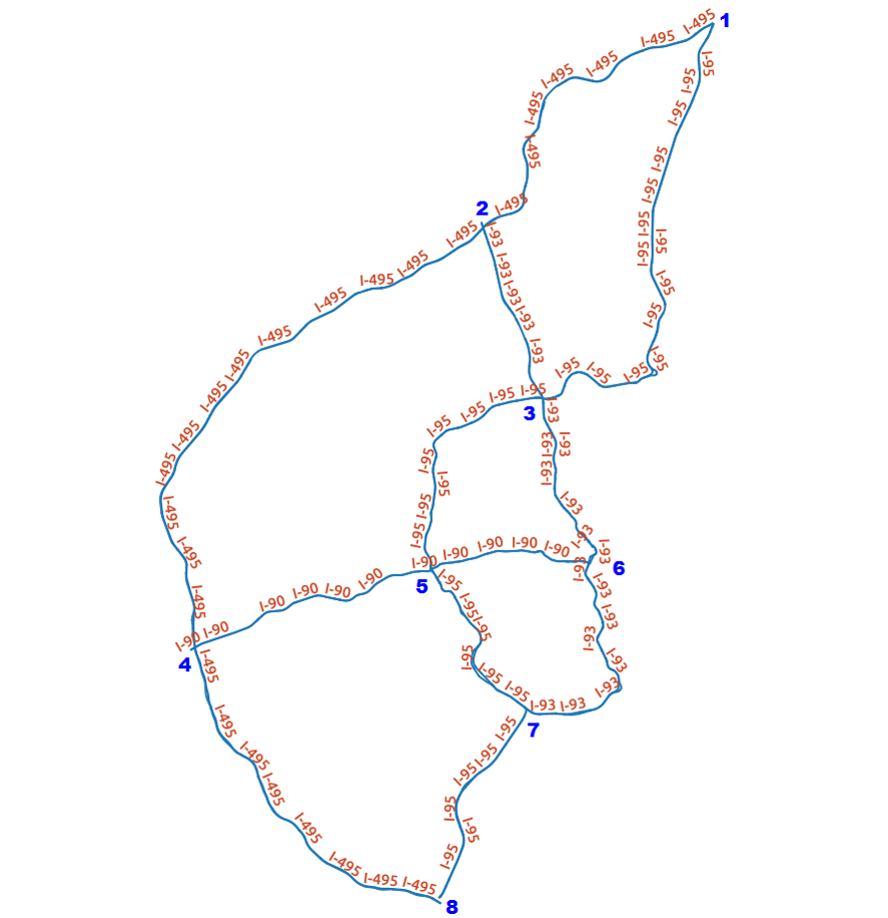}
\caption{}
\label{fig:EMA-small-subnet}
\end{subfigure}
\caption{ (a) Typical Traffic in a network from Google Maps; (b) Interstate highway subnetwork of EMA\cite{zhang_price_2016} }%
\label{fig:EMA Interstate}%
\end{figure}

Let us consider the traffic network as a directed graph (Fig.
\ref{fig:EMA Interstate}). Based on the traffic intensity on each link we can
categorize the links into three modes: low, medium, and high traffic links. We can then
assign different standard drive cycles to each link \cite{qiao_vehicle_2016} (HWFET $\rightarrow$ low traffic links, UDDS $\rightarrow$ medium traffic links, and NYC $\rightarrow$ high traffic links).
Qiao et al \cite{qiao_vehicle_2016}
modeled a PHEV20 (PHEV with 20 miles of AER) in \textit{PSAT} and calculated the average electrical energy ($\mu_{CD}$) and gas ($\mu_{CS}$) used to drive
one mile under CD and CS modes respectively under each of theses drive cycles
(Table \ref{tab: conversion factors}):%

\[%
\begin{array}
[c]{lr}%
\mu_{CD_{ij}}=\frac{d_{ij}}{P_{batt_{ij}}}\quad\text{,} & \mu_{CS_{ij}}%
=\frac{d_{ij}}{\dot{m}_{batt_{ij}}}\quad\text{,}%
\end{array}
\]
where $d_{ij}$ is the length of link $(i,j)$. By knowing $\mu_{CD_{ij}}$ and
$\mu_{CS_{ij}}$ on each link, as well as the network topology (length of each
link), we can determine the average fuel consumption rate ($\dot{m}%
_{batt_{ij}}$) and electrical power demand from the battery ($P_{batt_{ij}}$)
on each link $(i,j)$. We can then use (\ref{eqn:cost function}) to calculate the total energy cost for each trip.

\begin{table}[h]
\caption{Energy consumption of various drive cycles \cite{qiao_vehicle_2016}%
}%
\label{tab: conversion factors}
\centering
\resizebox{1\columnwidth}{!}{\renewcommand{\arraystretch}{0.9}
\begin{tabular}{llllll}
\hline
Vehicle Type & Symbol     & Unit     & HWFET & UDDS & NYC  \\ \hline
PHEV20       & $\mu_{CD}$ & $mi/kWh$ & 5.7   & 6.2  & 4.2  \\
& $\mu_{CS}$ & $mi/gal$ & 58.6  & 69.4 & 45.7 \\ \hline
\end{tabular}
}\end{table}



\section{Single Vehicle Eco-Routing}

\label{sec: single vehicle routing} 

In this section, we first review the CDF
eco-routing approach \cite{sun_save_2016,qiao_vehicle_2016}. Next, we propose the
CRPTC algorithm to
solve the eco-routing problem.

\subsection{Problem Formulation}

We model the traffic network as a directed graph $G=(\mathcal{N},\mathcal{A})$
with $\mathcal{N}={1,...,n}$ and $|\mathcal{A}|=m$ with the arc (link) connecting node \textit{i} to
\textit{j} denoted by $(i,j)\in\mathcal{A}$. The set of nodes that are
incoming/outgoing to node \textit{i} are defined as: $\mathcal{I}%
(i)=\{j\in\mathcal{N}|(j,i)\in\mathcal{A}\}$ and $\mathcal{O}(i)=\{j\in
\mathcal{N}|(i,j)\in\mathcal{A}\}$, respectively. We consider the 
single-origin-single-destination eco-routing problem where origin and
destination nodes are denoted by 1 and \textit{n} respectively. The energy
cost consumed by the vehicle on link $(i,j)$ is denoted by $c_{ij}$. We use
$E_{i}$ to represent the vehicle's residual battery energy at node $i$.
Moreover, we denote the selection of arc $(i,j)$ by $x_{ij}\in\{0,1\}$. The
problem objective is to determine a path from node 1 to \textit{n} so as to
minimize the total energy cost consumed by the vehicle to reach the
destination. We consider two approaches to solve this problem as follows.

\subsubsection{Charge Depleting First (CDF): A Hybrid-LP Relaxation Approach}

In this approach, we assume that the PHEV always starts every trip in  the CD
mode and uses electricity to drive the vehicle until it drains all the energy
out of the battery pack. Afterwards, it switches to the CS mode and starts
using gas to drive the vehicle. Even though this is not the most accurate
approach to solve the problem, it eliminates the need of using complicated
control strategies for a PHEV power-train to switch between Internal Combustion Engine (ICE) and electric
motors \cite{qiao_vehicle_2016}. As a result, we can formulate the eco-routing
problem using a mixed-integer nonlinear programming (MINLP) framework as follows:

%

\begin{gather}
\min_{x_{ij},i,j\in\mathcal{N}}\sum_{i=1}^{n}\sum_{j=1}^{n}c_{ij}%
x_{ij}\label{eqn: CDF objective}\\%
\label{eqn: CDF cost}
\begin{small}
\begin{array}
[c]{lc}%
s.t. & c_{ij}=%
\begin{cases}
C_{gas}\dfrac{d_{ij}}{\mu_{CS_{ij}}}; & E_{i}\leq0\\
C_{ele}\dfrac{d_{ij}}{\mu_{CD_{ij}}}; & E_{i}\geq\frac{d_{ij}}{\mu_{CD_{ij}}%
}\\
C_{ele}E_{i}+C_{gas}\dfrac{d_{ij}-\mu_{CD_{ij}}E_{i}}{\mu_{CS_{ij}}}; &
\text{otherwise}%
\end{cases}
\end{array}
\end{small}
\\%
\begin{array}
[c]{lc}%
E_{j}=\sum\limits_{i\in\mathcal{I}(j)}(E_{i}-\dfrac{d_{ij}}{\mu_{CD_{ij}}%
})x_{ij}, & \text{for }j=2,...,n,
\end{array}
\\%
\label{eqn: flow_CDF_1}
\begin{array}
[c]{cr}%
\sum\limits_{j\in\mathcal{O}(i)}x_{ij}-\sum\limits_{j\in\mathcal{I}(i)}%
x_{ij}=b_{i}, & \text{for each }i\in\mathcal{N}%
\end{array}
\\
\label{eqn: flow_CDF_2}
b_{1}=1,b_{n}=-1,b_{i}=0,\text{for }i\neq1,n\\
x_{ij}\in\{0,1\}
\end{gather}
where $E_{i}$ is the remaining electrical
energy at node \textit{i}, and $\mu_{CD_{ij}}$ and $\mu_{CS_{ij}}$ are the
conversion factors taken from Table \ref{tab: conversion factors}, which are
functions of the traffic intensity on each link $(i,j)$. Note that
(\ref{eqn: flow_CDF_1})-(\ref{eqn: flow_CDF_2}) are the flow conservation
constraints \cite{bertsimas_introduction_1997}. We assume that the vehicle has
enough gas and electrical power to complete the trip and that $E_{1}\geq0$.
Knowing the traffic density on each link, problem (\ref{eqn: CDF objective}) was
solved using Dijkstra's algorithm \cite{dijkstra_note_1959} in
\cite{qiao_vehicle_2016}. 

In what follows, we propose an alternative solution to this
problem which we call Hybrid-LP Relaxation. In this approach we reduce
the MINLP problem (\ref{eqn: CDF objective}) to a simpler problem which can be
solved using a combination of linear programming (LP) and a simple dynamic
programming-like algorithm, in order to guarantee global convergence. The
nonlinearities of the problem arise in (\ref{eqn: CDF cost}) where $c_{ij}$ is
a function of $x_{ij}$. We show that we can reduce this piecewise constant function to a
constant function, and the MINLP can be converted to a LP by using the properties
of the minimum cost flow problem \cite{hillier2012introduction}. The proposed
algorithm is as follows:


\begin{enumerate}
\item Find the shortest path (or shortest time path) and calculate the energy
cost on this path and set it to $\rho$. 

\item From the origin, construct all paths reaching node $p$ such that
$E_{p}\leq0$ and stop constructing the path at this node. Disregard the paths
with a total energy cost greater than $\rho$ and save the remaining paths in a matrix.

\item From (\ref{eqn: CDF cost}), the cost function for the paths outgoing
from node $p$ to $n$ in the previous step is given by:
\[
c_{ij}=C_{gas}\dfrac{d_{ij}}{\mu_{CS_{ij}}}%
\]

\item Assuming knowledge of traffic modes on each link, the least energy cost
path from node $p$ in step 2 to the destination node can be found from:
\begin{gather}
\min_{x_{ij},i,j\in\mathcal{N}}\sum_{i=1}^{n}\sum_{j=1}^{n}c_{ij}%
x_{ij}\label{eqn: hybrid LP}\\%
\begin{array}
[c]{lc}%
s.t. & c_{ij}=C_{gas}\dfrac{d_{ij}}{\mu_{CS_{ij}}}%
\end{array}
\\%
\begin{array}
[c]{cr}%
\sum\limits_{j\in\mathcal{O}(i)}x_{ij}-\sum\limits_{j\in\mathcal{I}(i)}%
x_{ij}=b_{i}, & \text{for each }i\in\mathcal{N}%
\end{array}
\\
\label{eqn: hybrid LP constraint}
b_{p}=1,b_{n}=-1,b_{i}=0,\text{for }i\neq p,n\\
x_{ij}\in\{0,1\}
\end{gather}
Note that constraint (\ref{eqn: hybrid LP constraint}) ensures that by solving (\ref{eqn: hybrid LP}), we are finding the optimal path from \(p\) to \(n\).

\item Using the property the minimum cost flow problem
\cite{hillier2012introduction}, problem (\ref{eqn: hybrid LP}) is equivalent to
an LP problem with the integer restriction of $x_{ij}$ relaxed:
\begin{equation}
0\leq x_{ij}\leq1
\end{equation}

\item Find the path from node 1 to $p$ with the least energy cost. By the
principle of optimality, the optimal path from 1 to $n$ is the one determined
in this manner followed by the path selected by steps 4 and 5 from node $p$ to $n$.

\item Find the paths in step 6 for all nodes $p$ such that $E_{p}\leq0$, then
choose the one with the minimum energy cost. The selected path would be the
minimum energy cost path.

\item If there are paths without any node such that $E_{p}\leq0$ (generated at
step 2), compare their cost function values with the cost functions in step 6.
The optimal route is the minimum among them.
\end{enumerate}
\begin{figure}[h]
\centering
\includegraphics[width=\columnwidth]{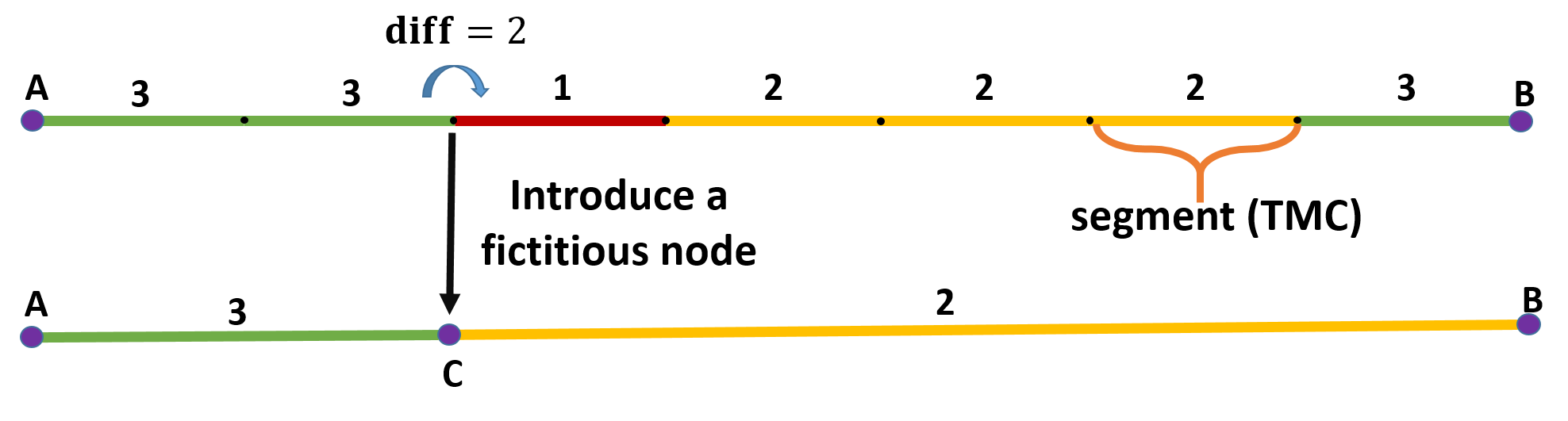}\caption{{Procedure
for introducing new fictitious nodes into a link based on the traffic modes on
the link's segments}}%
\label{fig:Link-Segment}%
\end{figure}

\subsubsection{Combined Routing and Power-Train Control (CRPTC)}
Based on Table \ref{tab: conversion factors}, the CD mode has the best
efficiency on medium traffic links. As such, if we always consider using
the CD mode at the beginning of each trip and then switch to the CS mode when we
run out of battery, we miss the opportunity to harness the effectiveness of the CD mode on medium traffic links towards the end of the route. With this
motivation, we propose a new algorithm which finds the routing decisions as
well as the PT controller decision to switch between CD and CS modes. Let
$y_{ij}\in\lbrack0,1]$ be an additional decision variable on link $(i,j)$
which represents the fraction of the link's length over which we use the CD
mode (thus, if we only use the CD mode over link $(i,j)$, then $y_{ij}=1$).
Considering the new decision variable, we can formulate the CRPTC problem as
follows:
\begin{gather}
\min_{x_{ij},y_{ij},i,j\in\mathcal{N}}\sum_{i=1}^{n}\sum_{j=1}^{n}%
[c_{gas}\dfrac{d_{ij}}{\mu_{CS_{ij}}}(1-y_{ij})+c_{ele}\dfrac{d_{ij}}%
{\mu_{CD_{ij}}}y_{ij}]x_{ij}\label{eqn:CRPTC Obj}\\%
\begin{array}
[c]{cr}%
s.t.\sum\limits_{j\in\mathcal{O}(i)}x_{ij}-\sum\limits_{j\in\mathcal{I}%
(i)}x_{ij}=b_{i}, & \text{for each }i\in\mathcal{N}%
\end{array}
\\
b_{1}=1,b_{n}=-1,b_{i}=0,\text{for }i\neq1,n\\
\label{eqn:CRPTC energy cons}
\sum_{i=1}^{n}\sum_{j=1}^{n}\dfrac{d_{ij}}{\mu_{CD_{ij}}}y_{ij}x_{ij}\leq
E_{1}\\
\begin{array}{lr}
x_{ij}\in\{0,1\},&
y_{ij}\in\lbrack0,1]
\end{array}
\end{gather}
Note that constraint (\ref{eqn:CRPTC energy cons}) ensures that the total
electrical energy used in the CD mode would be less than the initial available
energy in the battery ($E_{1}$). Since we have the term $x_{ij}y_{ij}$ in the
problem formulation, this is a MINLP problem and we may not be able to determine a global optimum. Hence, we
transform (\ref{eqn:CRPTC Obj}) into a mixed integer linear programming
problem (MILP) by introducing an intermediate decision variable $z_{ij}%
=x_{ij}y_{ij}$. We can then use the inequalities in (\ref{eqn:Zij inequalities}) to
transform the existing MINLP problem (\ref{eqn:CRPTC Obj}) into a MILP problem as follows:
\begin{gather}
\min_{x_{ij},y_{ij},z_{ij},i,j\in\mathcal{N}}\sum_{i=1}^{n}\sum_{j=1}%
^{n}(c_{gas}\dfrac{d_{ij}}{\mu_{CS_{ij}}}x_{ij}+(c_{ele}\dfrac{d_{ij}}%
{\mu_{CD_{ij}}}-c_{gas}\dfrac{d_{ij}}{\mu_{CS_{ij}}})z_{ij}%
)\label{eqn:Combined Single Vehicle energy}\\%
\begin{array}
[c]{cr}%
s.t.\sum\limits_{j\in\mathcal{O}(i)}x_{ij}-\sum\limits_{j\in\mathcal{I}%
(i)}x_{ij}=b_{i}, & \text{for each }i\in\mathcal{N}%
\end{array}
\\
b_{1}=1,b_{n}=-1,b_{i}=0,\text{for }i\neq1,n
\end{gather}
\begin{gather}
\sum_{i=1}^{n}\sum_{j=1}^{n}\dfrac{d_{ij}}{\mu_{CD_{ij}}}z_{ij}\leq E_{1}\\
\begin{array}{crc}
\label{eqn:Zij inequalities}
z_{ij}\geq 0 ,&    &
z_{ij}\leq y_{ij}\\
z_{ij}\leq x_{ij},&    &
z_{ij}\geq y_{ij}-(1-x_{ij})
\end{array}\\
\begin{array}{lr}
x_{ij}\in\{0,1\},&
y_{ij}\in\lbrack0,1]
\end{array}
\end{gather}
This is a MILP problem which can be solved to determine a global optimum.


\section{Numerical results}

\label{sec: numerical results} In order to evaluate the performance of the
proposed algorithm, we conduct a data-driven case study using the actual
traffic data from the EMA road network collected by
\textit{INRIX} \cite{zhang2018price,zhang2017data}. A sub-network including
the interstate highways of EMA (Fig. \ref{fig:EMA-small-subnet}) is chosen for
the case study. 
Details regarding this sub-network can be found in
\cite{zhang2018price,zhang2017data}. 

\subsection{Performance Measurement Baseline}

For measuring the 
performance of each algorithm, we need to define a
baseline against which to compare the energy cost obtained. We consider two
different baselines to compare with our energy-optimal algorithm: time optimal
paths and the actual paths from the INRIX dataset.

\subsubsection{Energy-Optimal Path vs. Time-Optimal Path}

In this approach, we first find the time-optimal paths using historical
traffic data as follows:
\begin{gather}
\min_{x_{ij},i,j\in\mathcal{N}}\sum_{i=1}^{n}\sum_{j=1}^{n}t_{ij}x_{ij}\\%
\begin{array}
[c]{lcc}%
s.t. & \label{eqn: travel time}t_{ij}=\dfrac{d_{ij}}{\bar{v}_{ij}}, & x_{ij}\in\{0,1\}%
\end{array}
\end{gather}
where $\bar{v}_{ij}$ and $d_{ij}$ are the average speed and length of link
$(i,j)$ respectively, and $t_{ij}$ is the travel time over link $(i,j)$. We
can then determine the energy costs for traveling through the shortest time
path using (\ref{eqn: CDF cost}) and compare the energy costs of traveling
through the optimal energy paths (under CDF and CRPTC) with the costs of
traveling through the minimum time path.

\subsubsection{Energy-Optimal Path vs. Actual Routes}

Using the INRIX dataset, Zhang et al \cite{zhang2018price,zhang2017data}
transformed the average speeds to vehicle flows on each link. They then solved
an inverse optimization problem \cite{bertsimas_data-driven_2015} to find the origin-destination (O-D) matrices of the network including the probabilities
of going through each route. They have reported these values based on four
different time periods: AM (6 am to 9 am), MD (9 am to 3pm), PM (3 pm to 6pm),
and NT (6 pm to 6 am). The values are calculated on an average sense over all
days in April 2012. An example is shown in Table \ref{tab: Actual Routes} for O-D pair \((1,5)\).
\begin{table}[h]
\caption{Actual routes and their probabilities \cite{zhang2018price}}%
\label{tab: Actual Routes}
\centering
\resizebox{1\columnwidth}{!}{\renewcommand{\arraystretch}{1}
\begin{tabular}{|c|c|c|c|c|c|c|}
\hline
Origin & Dest. & Route                                                                  & AM \% 		& MD\%		& PM\%		& NT\%		\\ \hline
1      & 5           & 1$\rightarrow$2$\rightarrow$3$\rightarrow$5                            & 19.7 		 & 18.7 	 & 21		 & 16.9		\\ \hline
1      & 5           & 1$\rightarrow$3$\rightarrow$5                                          & 65.3 		 & 73		 & 61		 & 73           \\ \hline
1      & 5           & 1$\rightarrow$2$\rightarrow$$\rightarrow$3$\rightarrow$6$\rightarrow$5 & 15 		 & 8.3		 & 18		 & 10.1             \\ \hline
\end{tabular}
}\end{table}

Using the actual routes, we calculate the expected energy costs and travel
times for each O-D pair based on the traffic information at any given time. We use the following equations to calculate the total expected
energy costs and traveling times for traveling through each O-D pair $(i,j)$:
\begin{equation}
\begin{array}{lr}
E(C_{ij}^{act})=\sum_{k=1}^{m}c_{ij}^{k}p_{ij}^{k}, & E(\tau_{ij}^{act})=\sum_{k=1}^{m}t_{ij}^{k}p_{ij}^{k}%
\end{array}
%
\\  
\end{equation}
where $E(\cdot)$ is the expected value, $C_{ij}^{act}$ and $\tau_{ij}^{act}$
are the total energy cost and traveling time for traveling from node $i$ to
node $j$ from (\ref{eqn: CDF cost}) and (\ref{eqn: travel time}), $m$ is the
number of possible routes, $c_{ij}^{k}$ and $t_{ij}^{k}$ are
the energy cost and travel time for going through $(i,j)$ using the $k$th possible
route, and $p_{ij}^{k}$ is the probability of the $k$th route for traveling
through link $(i,j)$.

\subsection{Data Preprocessing}

Each link consists of a number of road segments, and the average speed may differ
over consecutive road segments. Meanwhile, to solve problem
(\ref{eqn:Combined Single Vehicle energy}) we need to know the traffic mode on
each link. Since different segments on a link have different traffic modes
(low, medium, and heavy traffic), we need to come up with a strategy to assign
a unified link mode to each link. Our approach is as follows:

\begin{enumerate}
\item Categorize each link segment into 3 modes based on the average
speed $(V_{ave})$ of the segment: mode 1 ($V_{ave}<20$), mode 2 ($20\leq
V_{ave}\leq40$), mode 3 ($V_{ave}>40$).

\item If the change in the traffic mode of two consecutive segments is 2
, we introduce a fictitious node at that point
into our network graph (Fig.
\ref{fig:Link-Segment}). 
\item Calculate the average mode of the segments of each link, and report the
value as the traffic mode of that link.
\end{enumerate}

\begin{figure*}[ptb]
\centering
\begin{subfigure}{.4\textwidth}
\centering
\includegraphics[width=1\columnwidth]{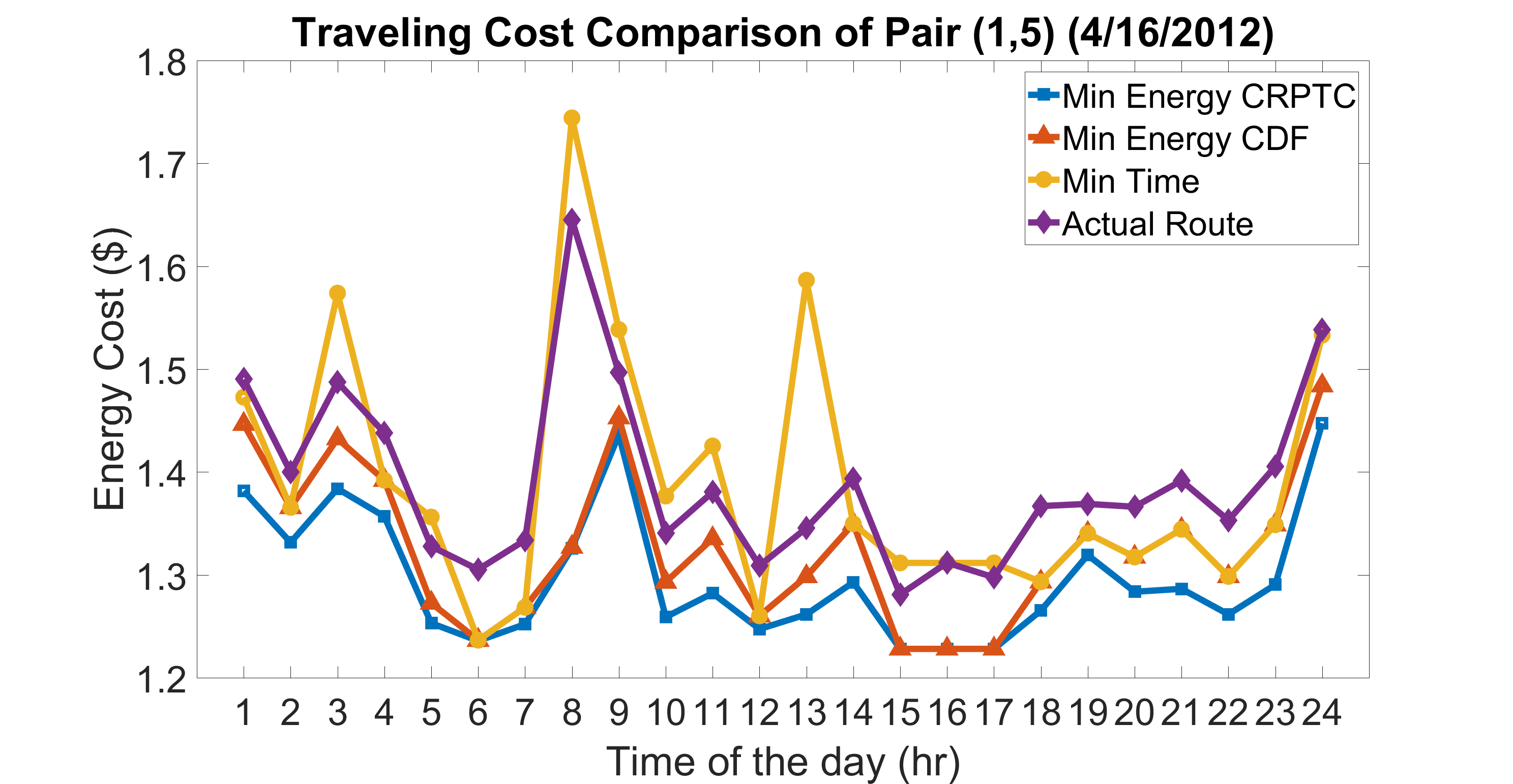}
\label{fig:Energy Cost Comparison Between All Routes (1,5) 4-17}
\end{subfigure}\begin{subfigure}{.4\textwidth}
\centering
\includegraphics[width=1\columnwidth]{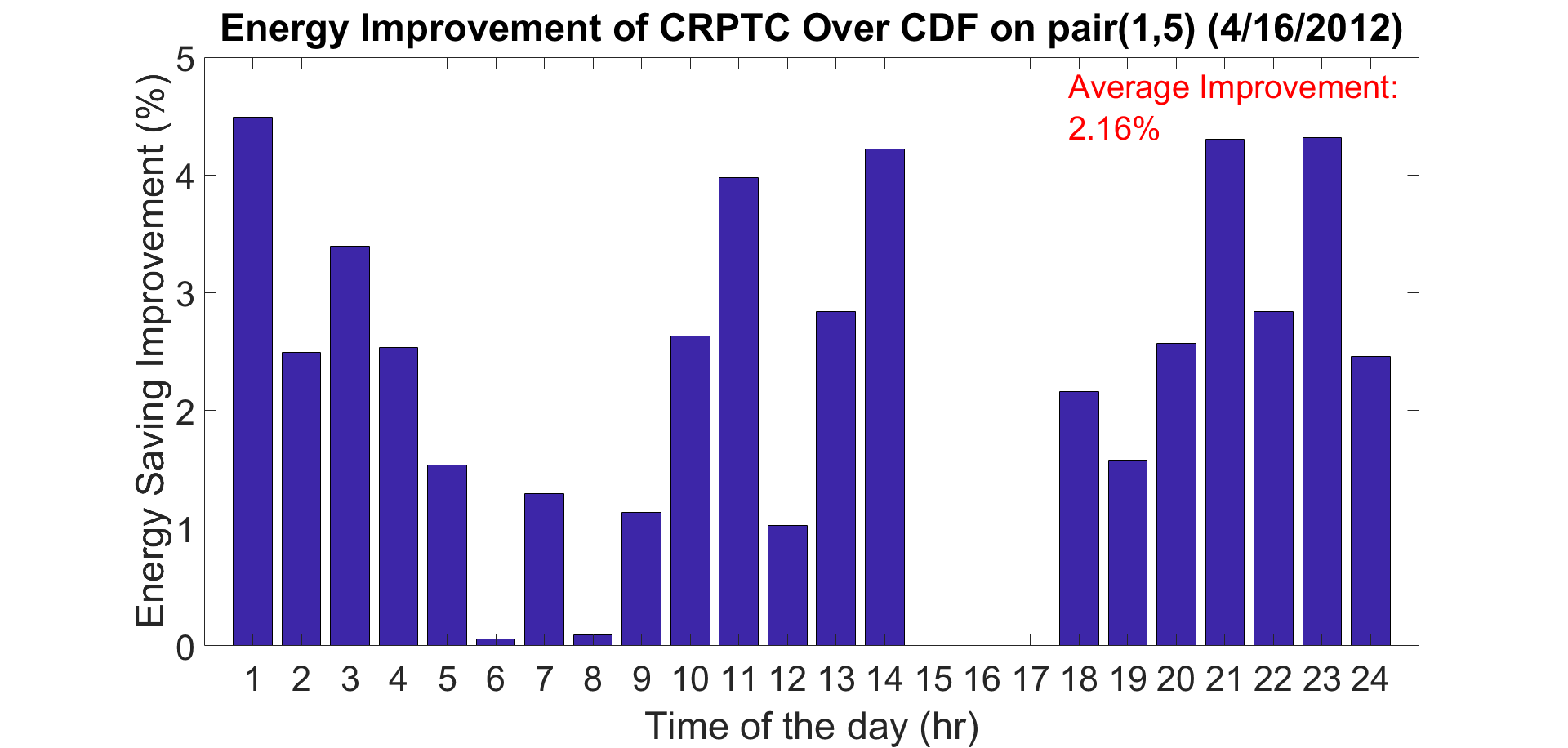}
\label{fig:Energy Cost of CRPTC VS CDF (1,5) 4-17}
\end{subfigure}
\begin{subfigure}{.4\textwidth}
\centering
\includegraphics[width=1\columnwidth]{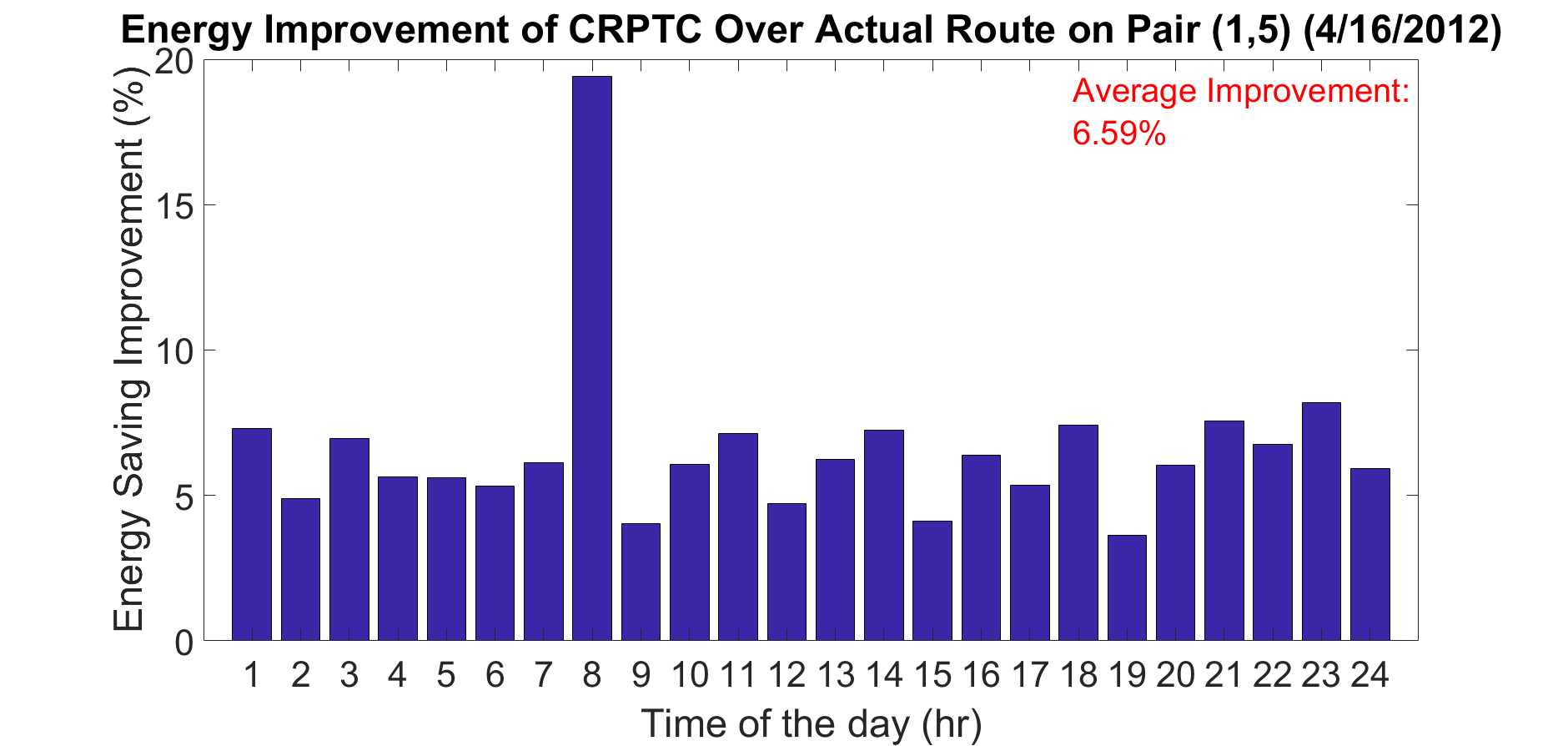}
\label{fig:Energy Cost CRPTC VS Actual Route (1,5) 4-17}
\end{subfigure}\begin{subfigure}{.4\textwidth}
\centering
\includegraphics[width=1\columnwidth]{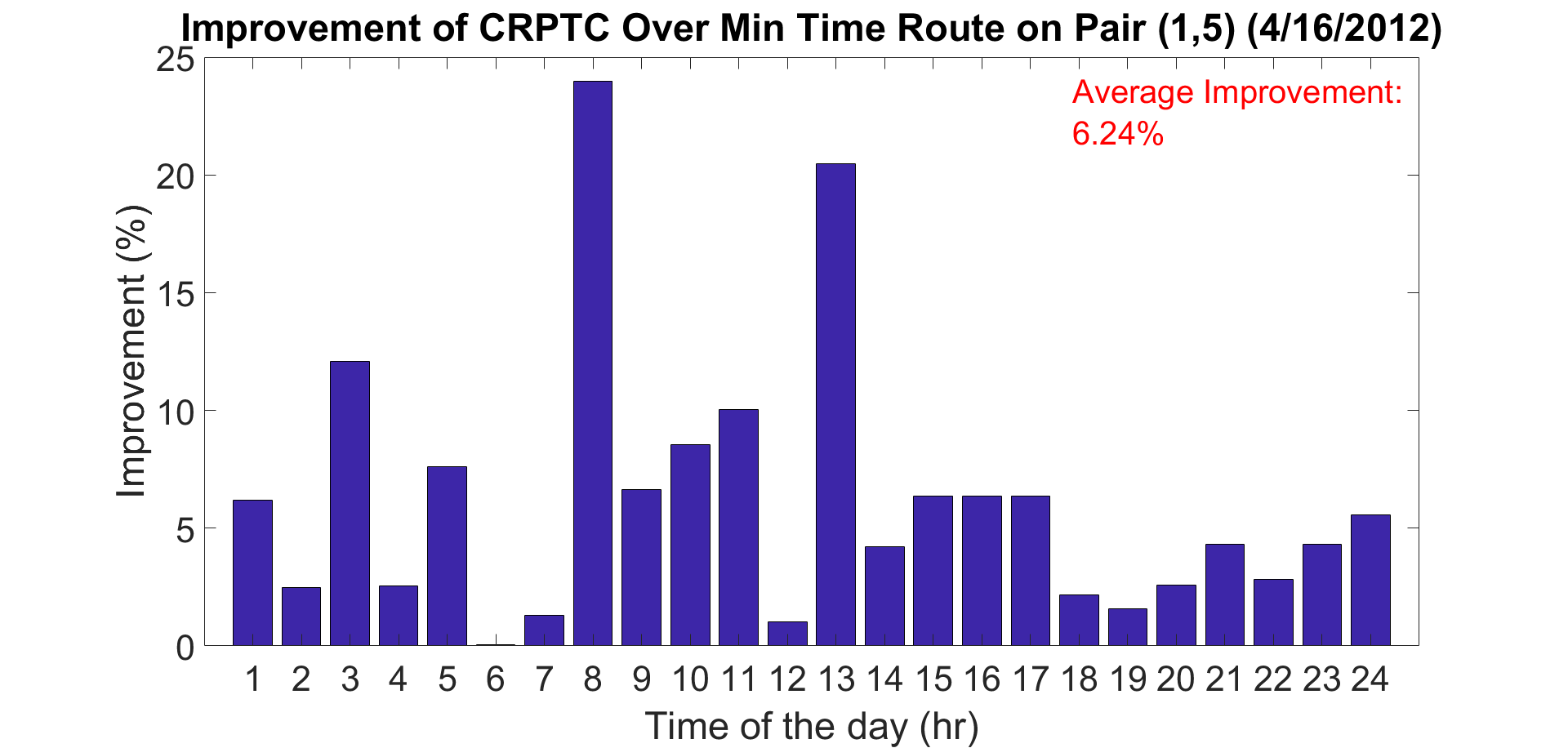}
\label{fig:Energy Cost CRPTC VS Min Time Route (1,5) 4-17}
\end{subfigure}
\par
\caption{Energy comparison plots for traveling from node 1 to 5 (04/16/2012)}%
\label{fig: Energy Comparison plots}%
\end{figure*}

\begin{figure*}[ptb]
\centering
\begin{subfigure}{.4\textwidth}
\centering
\includegraphics[width=1\columnwidth]{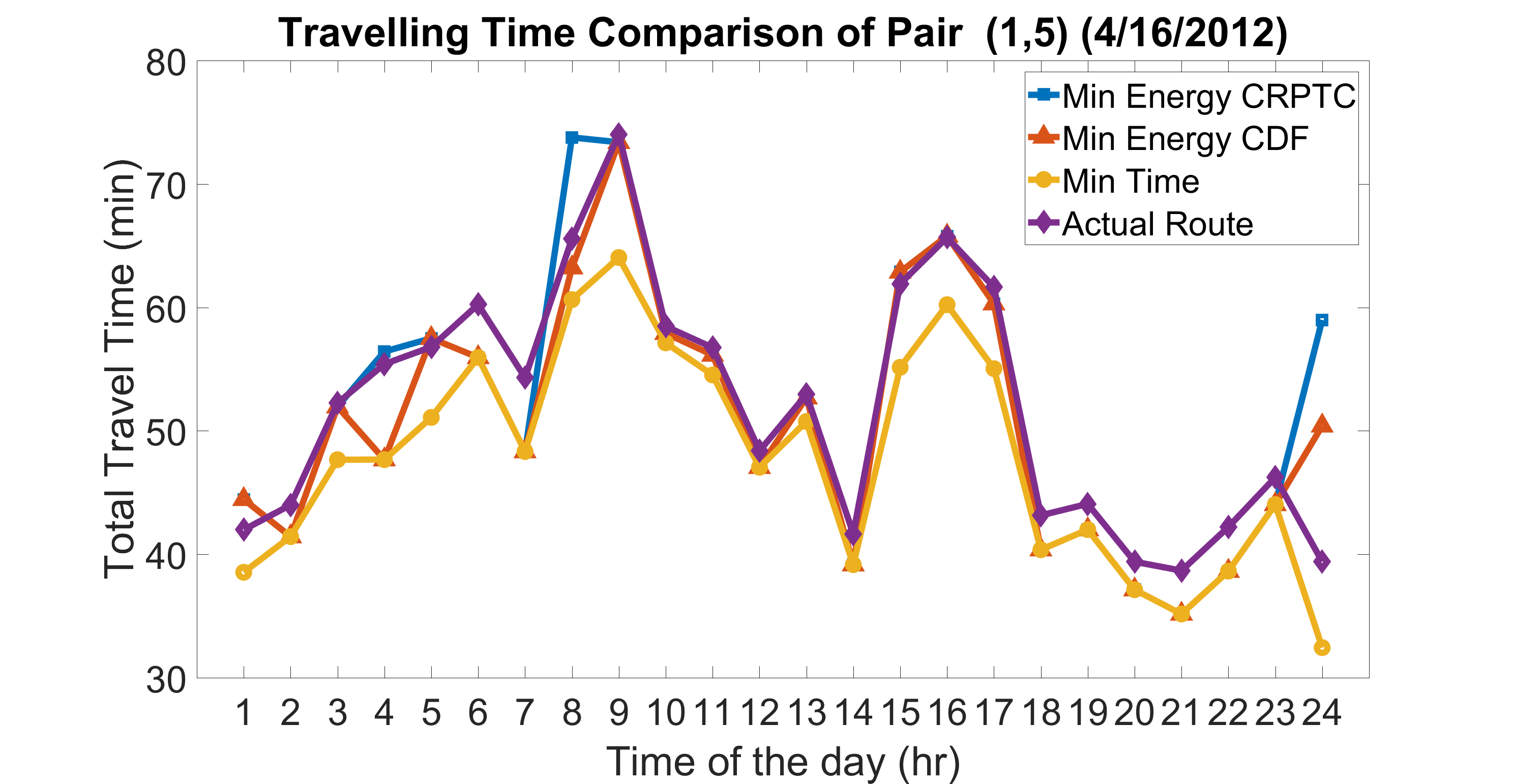}
\label{fig:Time Comparison Between All Routes (1,5) 4-16}
\end{subfigure}\begin{subfigure}{.4\textwidth}
\centering
\includegraphics[width=1\columnwidth]{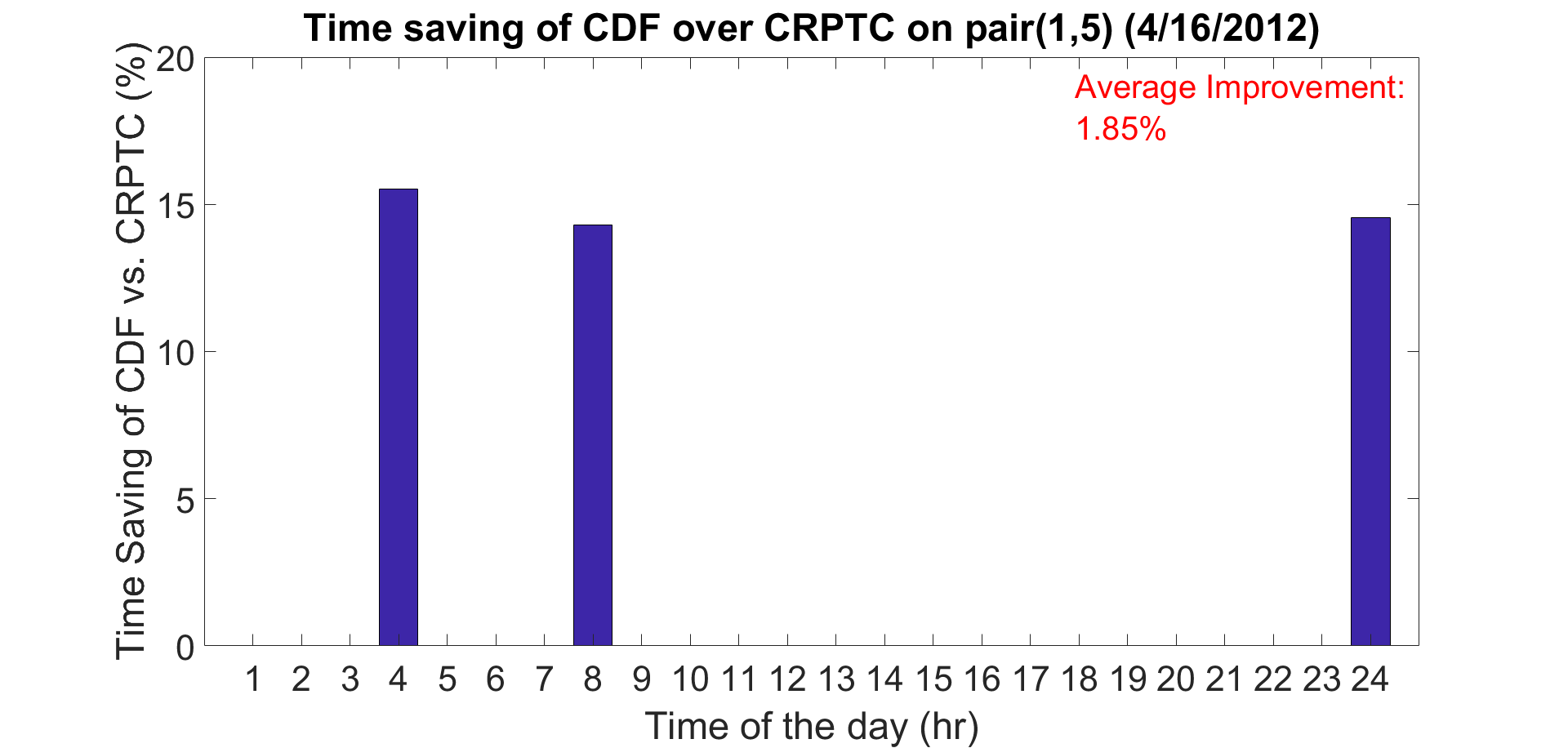}
\label{fig:Traveling Time of CRPTC VS CDF (1,5) 4-16}
\end{subfigure}
\begin{subfigure}{.4\textwidth}
\centering
\includegraphics[width=1\columnwidth]{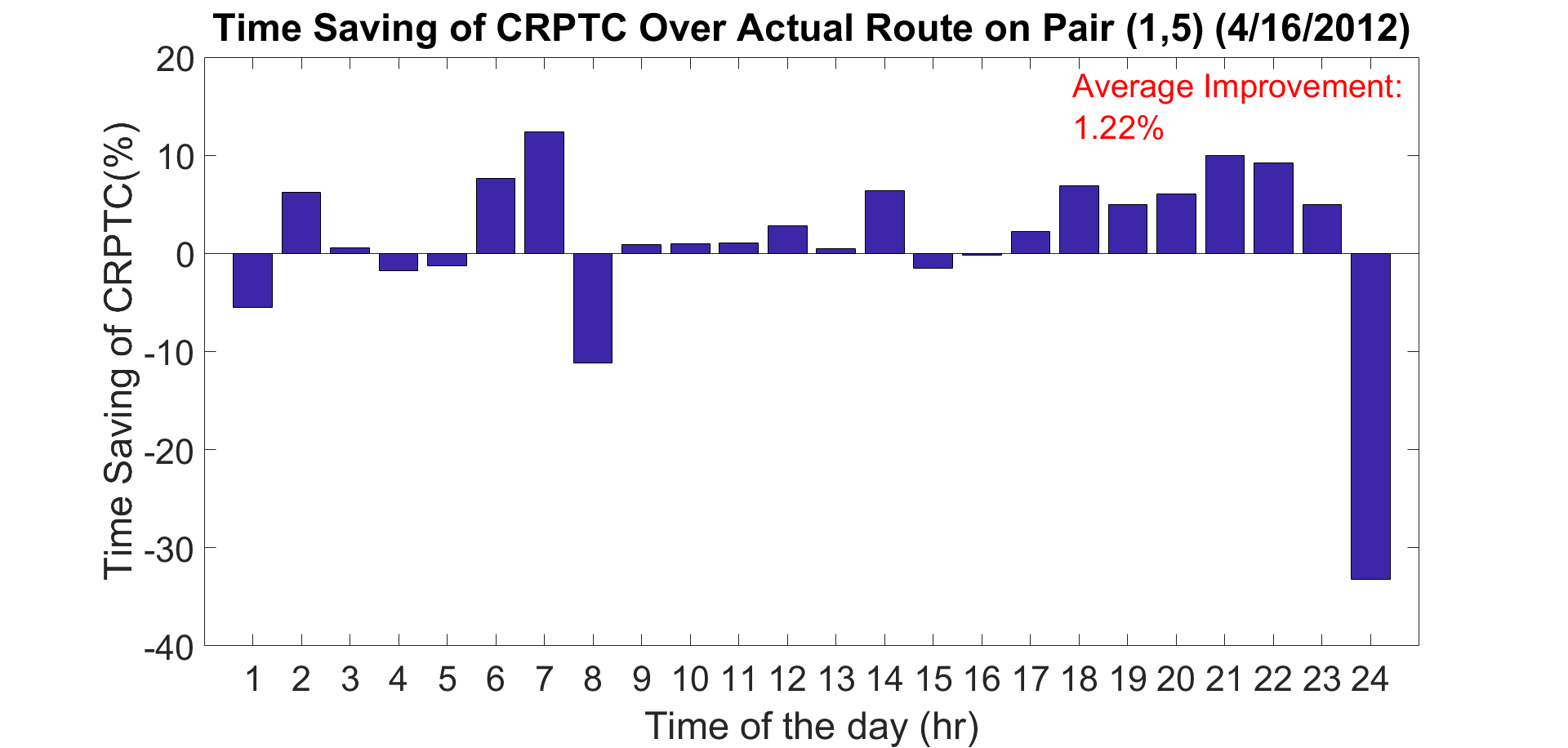}
\label{fig:Traveling Time of CRPTC VS Actual Route (1,5) 4-16}
\end{subfigure}\begin{subfigure}{.4\textwidth}
\centering
\includegraphics[width=1\columnwidth]{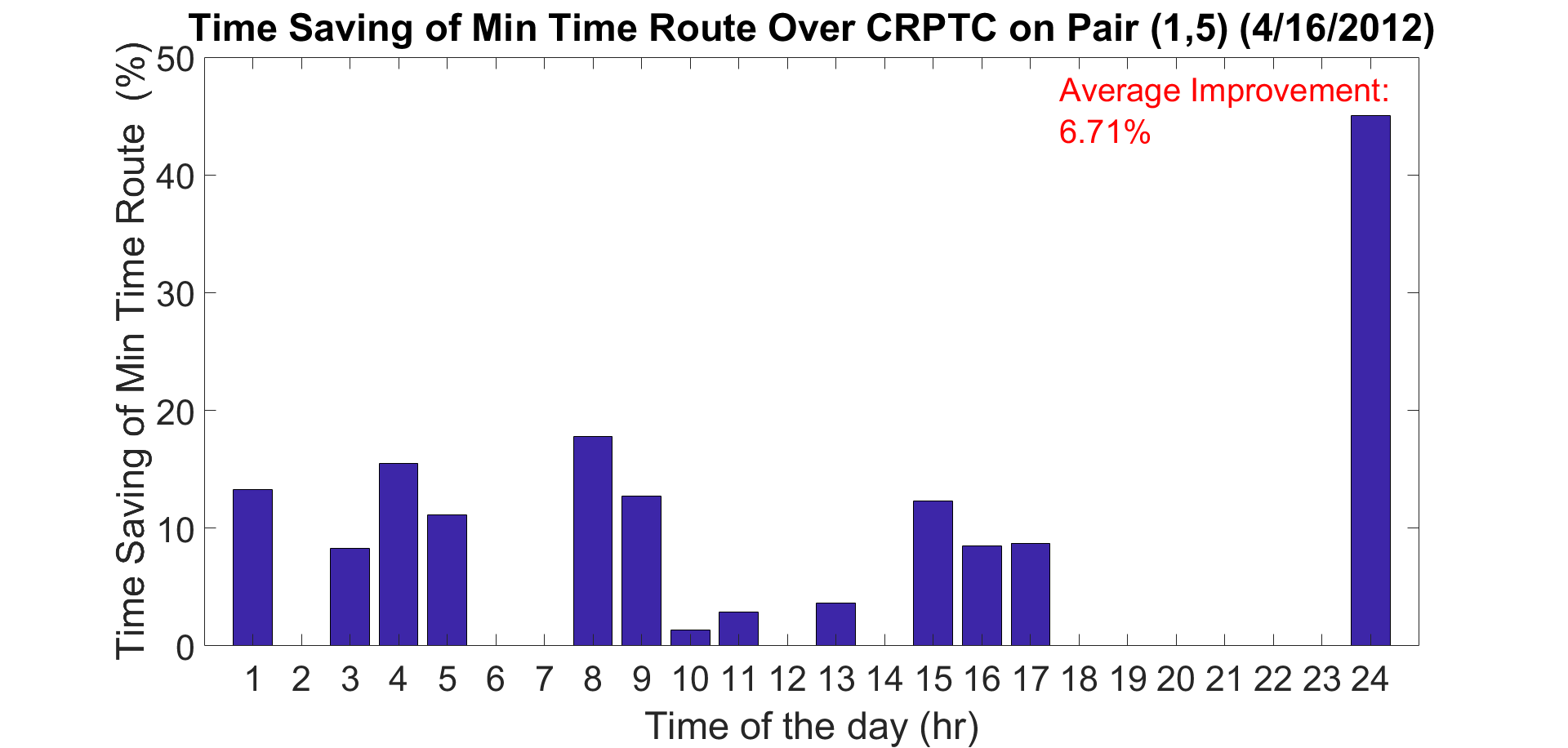}
\label{fig:Traveling Time of CRPTC VS Min Time Route (1,5) 4-16}
\end{subfigure}
\par
\caption{Traveling Time plots for traveling from node 1 to 5 (04/16/2012)}%
\label{fig: time Comparison plots}%
\end{figure*}
Considering this approach, we end up with a new adjacency matrix including
more nodes and edges than the original network. Using the new adjacency
matrix, we can solve problems (\ref{eqn: CDF objective}) and
(\ref{eqn:Combined Single Vehicle energy}).

\subsection{Comparison Results}

Using the INRIX dataset, optimal energy paths (CDF and CRPTC), optimal time paths, and
actual paths have been evaluated. We have also calculated the energy
consumption costs and travel times for each of these paths. The cost of
gas and electricity used in this case study are $C_{gas}=2.75\$/gal$ and
$C_{ele}=0.114\$/kWh$ respectively \cite{karabasoglu_influence_2013}.
Moreover, as in \cite{qiao_vehicle_2016}, we assumed the initial
available battery energy to be $E_{ini}=5.57kWh$. Energy and time
comparison plots for traveling thorough node 1 to 5 can be found in Figs.
\ref{fig: Energy Comparison plots} and \ref{fig: time Comparison plots}.

As expected, CRPTC performs better than CDF in terms of energy saving. We
repeated the same analysis for traveling from node 1 to 5 for a week in April
2012 using INRIX dataset and found that on average we can save 2.54\% in
terms of energy using CRPTC instead of CDF. We have also found that we can
save on average 5.02\% in terms of energy using CRPTC compared to the cost for
traveling through the min time route, and also 6.89\% compared to the cost for
traveling through the actual routes. The trade-off between energy saving and
travel time is such that on average it takes 4.79\% more time to travel from
node 1 to 5 while using the CRPTC route compared to traveling through the
minimum time route. It also takes 2.15\% longer using the CRPTC route compared
to the CDF route. We observe that users were not taking the time-optimal route
in 2012; this may be because in 2012 navigation systems with real time traffic
data were not as accessible to the public as they are today. In fact, we can
save both in terms of energy and time if we take the CRPTC route instead of
the actual routes users traveled at the time, and the average time saving is 2.15\%.


\subsection{Traffic Simulation}

Since we did not want to rely solely upon the historical traffic data to
validate our routing algorithm, we decided to simulate the traffic of the EMA
sub-network (Fig. \ref{fig:EMA Interstate}) using SUMO (Simulation of Urban
MObility) \cite{krajzewicz_recent_2012}. The flow data on each link are needed
to start a simulation in SUMO. Hence, we used the INRIX data to extract the
flow data for the subnetwork. The details of
theses calculations may be found in
\cite{zhang2018price,pourazarm_optimal_2016}.

\begin{figure}[h]
\centering
\begin{subfigure}{.39\textwidth}
\centering
\includegraphics[width=1\linewidth]{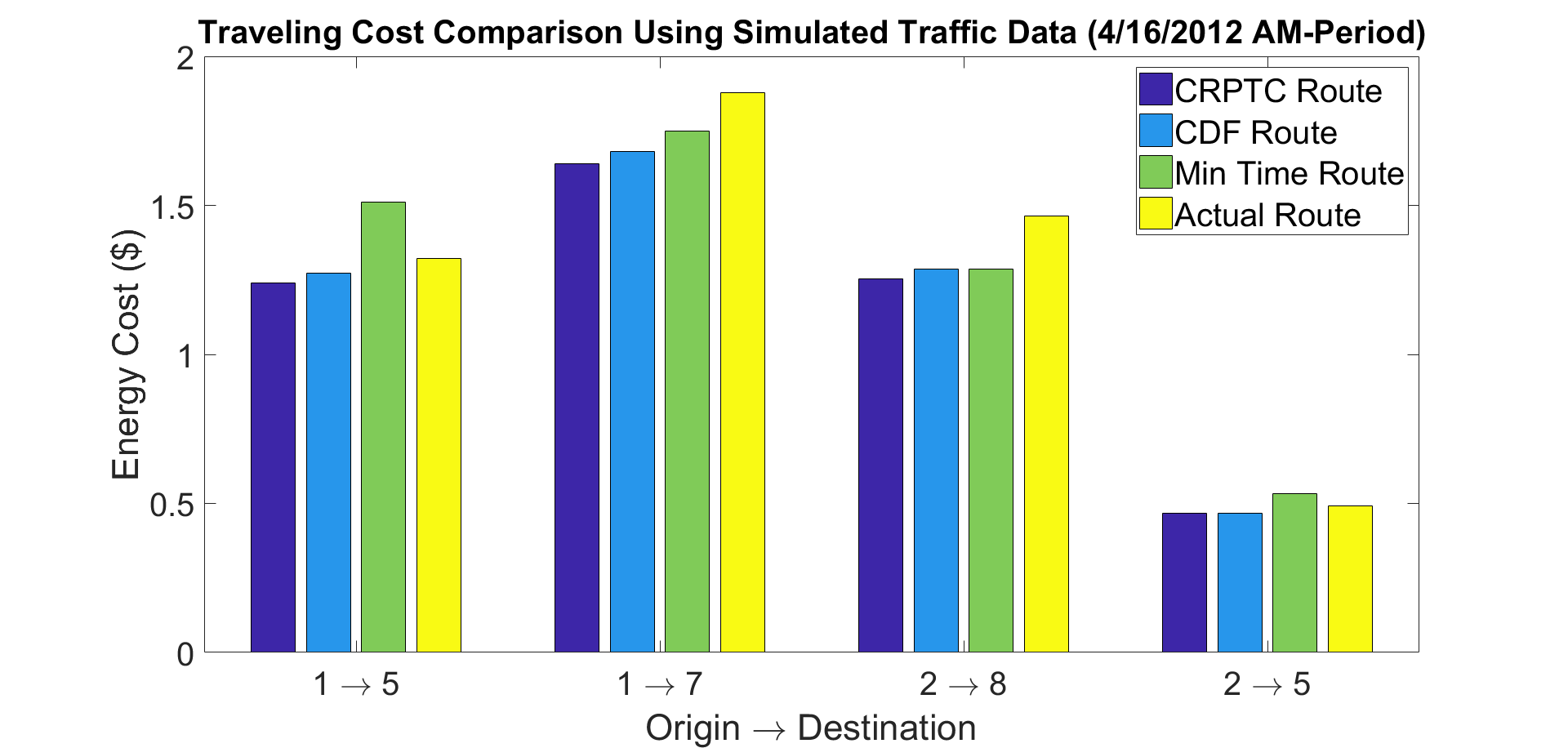}
\caption{Energy comparison plots for different O-D pairs}
\label{fig:SUMO Energy}
\end{subfigure}
\begin{subfigure}{.39\textwidth}
\centering
\includegraphics[width=01\linewidth]{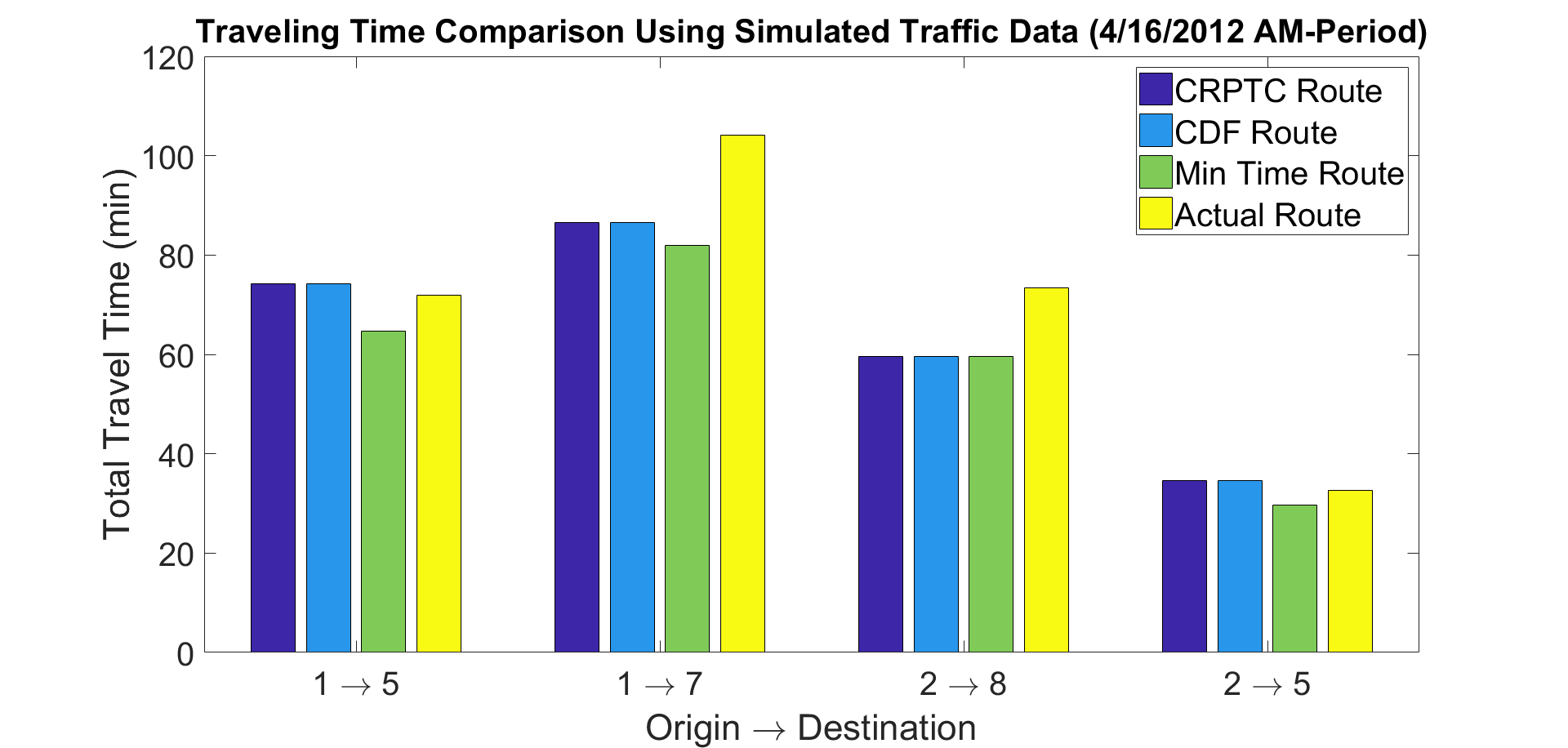}
\caption{ Time comparison plots for different O-D pairs}
\label{fig:SUMO Time}
\end{subfigure}
\caption{ Comparison plots using SUMO simulated traffic data}%
\end{figure}

Using extracted flows, we simulated traffic in SUMO. We then
aggregated every 5 segments in the map into a single link and recorded the
average speed of that new link. In this respect, we ended up having a graph
with 281 nodes and 300 edges. Using the new graph, we found the minimum time
route as well as the CRPTC and CDF energy-optimal routes, and calculated the
energy costs and traveling time for each of these routes for different O-D
pairs. Considering Table \ref{tab: Actual Routes}, we have also calculated the
expected energy costs and travel times over the actual paths taken
by drivers. The comparison results for April 16, 2012 between energy costs and
traveling time of 4 different O-D pairs are shown in Figs.
\ref{fig:SUMO Energy} and \ref{fig:SUMO Time}.



\section{Conclusions and Future Work}

\label{sec: conclusions}

In this paper, we reviewed the current approaches to solve the eco-routing
problem for PHEVs, and proposed a method to solve the minimum energy cost
problem for a single vehicle routing problem. The proposed CRPTC method is
capable of finding both the optimal path and the optimal switching strategy
between CD and CS modes on each link. Historical traffic data, as well as SUMO
simulations, were used to validate the performance of our algorithm. Numerical results show improvements in terms of the total energy cost using the CRPTC approach compared to the previous CDF method. We have also shown
that there is a trade-off between energy saving and time saving.

So far, we have not considered dynamically updating routing decisions at network nodes to account for sudden changes in traffic conditions (e.g., due to accidents). In ongoing
work, we are implementing and investigating dynamic eco-routing as well.
Moreover, we have so far solved the problem for a single vehicle scenario
with a known origin and destination. As a next step, we will consider
connectivity among vehicles and determine the social optimum for the network
considering a 100\% penetration rate of connected automated vehicle. Moreover, we plan to include multiple
vehicle architectures with different fuel consumption models and different
initial energies to the problem, as well as add charging stations in the
network to let vehicles recharge their batteries if necessary.

\bibliographystyle{IEEEtran}
\begin{tiny}
\bibliography{NEXTCAR}

\begin{thebibliography}{10}
\providecommand{\url}[1]{#1}
\csname url@samestyle\endcsname
\providecommand{\newblock}{\relax}
\providecommand{\bibinfo}[2]{#2}
\providecommand{\BIBentrySTDinterwordspacing}{\spaceskip=0pt\relax}
\providecommand{\BIBentryALTinterwordstretchfactor}{4}
\providecommand{\BIBentryALTinterwordspacing}{\spaceskip=\fontdimen2\font plus
\BIBentryALTinterwordstretchfactor\fontdimen3\font minus
  \fontdimen4\font\relax}
\providecommand{\BIBforeignlanguage}[2]{{%
\expandafter\ifx\csname l@#1\endcsname\relax
\typeout{** WARNING: IEEEtran.bst: No hyphenation pattern has been}%
\typeout{** loaded for the language `#1'. Using the pattern for}%
\typeout{** the default language instead.}%
\else
\language=\csname l@#1\endcsname
\fi
#2}}
\providecommand{\BIBdecl}{\relax}
\BIBdecl

\bibitem{bertsekas_dynamic_1995}
D.~P. Bertsekas, \emph{Dynamic programming and optimal control}.\hskip 1em plus
  0.5em minus 0.4em\relax Athena scientific Belmont, MA, 1995, vol.~1, no.~2.

\bibitem{braekers_vehicle_2016}
K.~Braekers, K.~Ramaekers, and I.~Van~Nieuwenhuyse,
  ``\BIBforeignlanguage{en}{The vehicle routing problem: {State} of the art
  classification and review},'' \emph{\BIBforeignlanguage{en}{Computers \&
  Industrial Engineering}}, vol.~99, pp. 300--313, Sep. 2016.

\bibitem{toth_vehicle_2002}
P.~Toth and D.~Vigo, \emph{The vehicle routing problem}.\hskip 1em plus 0.5em
  minus 0.4em\relax SIAM, 2002.

\bibitem{barth_environmentally-friendly_2007}
M.~Barth, K.~Boriboonsomsin, and A.~Vu, ``Environmentally-friendly
  navigation,'' in \emph{Intelligent {Transportation} {Systems} {Conference},
  2007. {ITSC} 2007. {IEEE}}.\hskip 1em plus 0.5em minus 0.4em\relax IEEE,
  2007, pp. 684--689.

\bibitem{boriboonsomsin_eco-routing_2012}
K.~Boriboonsomsin, M.~J. Barth, W.~Zhu, and A.~Vu, ``Eco-{Routing} {Navigation}
  {System} {Based} on {Multisource} {Historical} and {Real}-{Time} {Traffic}
  {Information},'' \emph{IEEE Transactions on Intelligent Transportation
  Systems}, vol.~13, no.~4, pp. 1694--1704, Dec. 2012.

\bibitem{andersen_ecotour:_2013}
O.~Andersen, C.~S. Jensen, K.~Torp, and B.~Yang, ``{EcoTour}: {Reducing} the
  {Environmental} {Footprint} of {Vehicles} {Using} {Eco}-routes,'' in
  \emph{2013 {IEEE} 14th {International} {Conference} on {Mobile} {Data}
  {Management}}, vol.~1, Jun. 2013, pp. 338--340.

\bibitem{yao_study_2013}
E.~Yao and Y.~Song, ``Study on eco-route planning algorithm and environmental
  impact assessment,'' \emph{Journal of Intelligent Transportation Systems},
  vol.~17, no.~1, pp. 42--53, 2013.

\bibitem{yang_stochastic_2014}
B.~Yang, C.~Guo, C.~S. Jensen, M.~Kaul, and S.~Shang, ``Stochastic skyline
  route planning under time-varying uncertainty,'' in \emph{2014 {IEEE} 30th
  {International} {Conference} on {Data} {Engineering}}, Mar. 2014, pp.
  136--147.

\bibitem{kubicka_performance_2016}
M.~Kubička, J.~Klusáček, A.~Sciarretta, A.~Cela, H.~Mounier, L.~Thibault,
  and S.-I. Niculescu, ``Performance of current eco-routing methods,'' in
  \emph{Intelligent {Vehicles} {Symposium} ({IV}), 2016 {IEEE}}.\hskip 1em plus
  0.5em minus 0.4em\relax IEEE, 2016, pp. 472--477.

\bibitem{guanetti_control_2018}
J.~Guanetti, Y.~Kim, and F.~Borrelli, ``Control of {Connected} and {Automated}
  {Vehicles}: {State} of the {Art} and {Future} {Challenges},''
  \emph{arXiv:1804.03757 [cs]}, Apr. 2018, arXiv: 1804.03757.

\bibitem{cela_energy_2014}
A.~Cela, T.~Jurik, R.~Hamouche, R.~Natowicz, A.~Reama, S.~I. Niculescu, and
  J.~Julien, ``Energy {Optimal} {Real}-{Time} {Navigation} {System},''
  \emph{IEEE Intelligent Transportation Systems Magazine}, vol.~6, no.~3, pp.
  66--79, 2014.

\bibitem{sun_save_2016}
Z.~Sun and X.~Zhou, ``To save money or to save time: {Intelligent} routing
  design for plug-in hybrid electric vehicle,'' \emph{Transportation Research
  Part D: Transport and Environment}, vol.~43, pp. 238--250, Mar. 2016.

\bibitem{qiao_vehicle_2016}
Z.~Qiao and O.~Karabasoglu, ``Vehicle {Powertrain} {Connected} {Route}
  {Optimization} for {Conventional}, {Hybrid} and {Plug}-in {Electric}
  {Vehicles},'' \emph{arXiv:1612.01243 [cs]}, Dec. 2016, arXiv: 1612.01243.

\bibitem{kamal_model_2013}
M.~A.~S. Kamal, M.~Mukai, J.~Murata, and T.~Kawabe, ``Model predictive control
  of vehicles on urban roads for improved fuel economy,'' \emph{IEEE
  Transactions on control systems technology}, vol.~21, no.~3, pp. 831--841,
  2013.

\bibitem{karabasoglu_influence_2013}
O.~Karabasoglu and J.~Michalek, ``Influence of driving patterns on life cycle
  cost and emissions of hybrid and plug-in electric vehicle powertrains,''
  \emph{Energy Policy}, vol.~60, pp. 445--461, 2013.

\bibitem{zhang_price_2016}
J.~Zhang, S.~Pourazarm, C.~G. Cassandras, and I.~C. Paschalidis, ``The price of
  anarchy in transportation networks by estimating user cost functions from
  actual traffic data,'' in \emph{2016 {IEEE} 55th {Conference} on {Decision}
  and {Control} ({CDC})}, Dec. 2016, pp. 789--794.

\bibitem{bertsimas_introduction_1997}
D.~Bertsimas, J.~N. Tsitsiklis, and J.~Tsitsiklis,
  \emph{\BIBforeignlanguage{English}{Introduction to {Linear} {Optimization}}},
  third printing edition~ed.\hskip 1em plus 0.5em minus 0.4em\relax Belmont,
  Mass: Athena Scientific, Feb. 1997.

\bibitem{dijkstra_note_1959}
E.~W. Dijkstra, ``\BIBforeignlanguage{en}{A note on two problems in connexion
  with graphs},'' \emph{\BIBforeignlanguage{en}{Numer. Math.}}, vol.~1, no.~1,
  pp. 269--271, Dec. 1959.

\bibitem{hillier2012introduction}
F.~S. Hillier, \emph{Introduction to operations research}.\hskip 1em plus 0.5em
  minus 0.4em\relax Tata McGraw-Hill Education, 2012.

\bibitem{zhang2018price}
J.~Zhang, S.~Pourazarm, C.~G. Cassandras, and I.~C. Paschalidis, ``The {Price}
  of {Anarchy} in {Transportation} {Networks}: {Data}-{Driven} {Evaluation} and
  {Reduction} {Strategies},'' \emph{Proceedings of the IEEE}, vol. 106, no.~4,
  pp. 538--553, Apr. 2018.

\bibitem{zhang2017data}
------, ``Data-driven estimation of origin-destination demand and user cost
  functions for the optimization of transportation networks,''
  \emph{IFAC-PapersOnLine}, vol.~50, no.~1, pp. 9680--9685, 2017.

\bibitem{bertsimas_data-driven_2015}
D.~Bertsimas, V.~Gupta, and I.~C. Paschalidis,
  ``\BIBforeignlanguage{en}{Data-driven estimation in equilibrium using inverse
  optimization},'' \emph{\BIBforeignlanguage{en}{Math. Program.}}, vol. 153,
  no.~2, pp. 595--633, Nov. 2015.

\bibitem{krajzewicz_recent_2012}
D.~Krajzewicz, J.~Erdmann, M.~Behrisch, and L.~Bieker, ``Recent development and
  applications of {SUMO}-{Simulation} of {Urban} {MObility},''
  \emph{International Journal On Advances in Systems and Measurements}, vol.~5,
  no. 3\&4, 2012.

\bibitem{pourazarm_optimal_2016}
S.~Pourazarm, C.~G. Cassandras, and T.~Wang, ``\BIBforeignlanguage{en}{Optimal
  routing and charging of energy-limited vehicles in traffic networks},''
  \emph{\BIBforeignlanguage{en}{Int. J. Robust. Nonlinear Control}}, vol.~26,
  no.~6, pp. 1325--1350, Apr. 2016.

\end{thebibliography}
\end{tiny}

\end{document}